\documentclass{amsart}
\usepackage{amssymb}
\theoremstyle{plain}
 \newtheorem{thm}{Theorem}
 \newtheorem{cor}{Corollary}
 \newtheorem{lem}{Lemma}
 
\theoremstyle{definition}

 \newtheorem{exmp}{Example}
\theoremstyle{remark}

 \newtheorem{ack}{Acknowledgment\ignorespaces}
 
\allowdisplaybreaks
\newcommand{\NaturalNumber}{\mathbb N}
\newcommand{\Integer}{\mathbb Z}
\newcommand{\RationalNumber}{\mathbb Q}
\newcommand{\RealNumber}{\mathbb R}
\renewcommand{\labelenumi}{(\roman{enumi})}

\begin{document}
\title[$n$-parameter semigroup]
{The set of common fixed points of an $n$-parameter continuous semigroup
 of mappings}
\author[T. Suzuki]{Tomonari Suzuki}
\date{}
\address{Department of Mathematics,
Kyu\-shu Institute of Technology,
1-1, Sen\-sui\-cho, To\-bata\-ku, Kita\-kyu\-shu 804-8550, Japan}
\email{suzuki-t@mns.kyutech.ac.jp}

\begin{abstract}
In this paper,
 using Kronecker's theorem,
 we discuss
 the set of common fixed points of
 an $n$-parameter continuous semigroup $\{ T(p) : p \in \RealNumber_+^n \}$
 of mappings.
We also discuss convergence theorems
 to a common fixed point of
 an $n$-parameter nonexpansive semigroup $\{ T(p) : p \in \RealNumber_+^n \}$.
\end{abstract}
\maketitle

\section{Introduction}
\label{SC:introduction}

Throughout this paper,
 we denote by
 $\NaturalNumber$,
 $\Integer$,
 $\RationalNumber$ and
 $\RealNumber$
 the sets of
 all positive integers,
 all integers,
 all rational numbers and
 all real numbers, respectively.
We put $\RealNumber_+^n = [0,\infty)^n$ and
 $$ e_j = (0,0,\cdots,0,0,\stackrel{(j)}{\mathstrut 1},0,0,\cdots,0)
 \in \RealNumber^n $$
 for $j \in \NaturalNumber$ with $1 \leq j \leq n$.

Let $C$ be a subset of a Banach space $E$, and
 let $T$ be a nonexpansive mapping on $C$, i.e.,
 $\| Tx - Ty \| \leq \| x - y \|$
 for all $x, y \in C$.
We know that $T$ has a fixed point
 in the case that $E$ is uniformly convex and
 $C$ is bounded, closed and convex;
 see \cite{REF:Browder1965_ProcNAS_3,REF:Gohde1965}.
See also
 \cite{REF:Baillon1979_1,
 REF:Browder1965_ProcNAS_4,REF:Kirk1965_AMMonth}
 and others.
We denote by $F(T)$ the set of fixed points of $T$.

Let $\tau$ be a Hausdorff topology on $E$.
A family of mappings $\{ T(p) : p \in \RealNumber_+^n \}$ is called
 an $n$-parameter $\tau$-continuous semigroup of mappings on $C$
 if the following are satisfied:
 \begin{enumerate}
 \renewcommand{\labelenumi}{(sg \arabic{enumi})}
 \item
 $ T(p+q) = T(p) \circ T(q) $ for all $p, q \in \RealNumber_+^n$;
 \item
 for each $x \in C$,
 the mapping $p \mapsto T(p)x $ from $\RealNumber_+^n$ into $C$ is continuous
 with respect to $\tau$.
 \end{enumerate}
As a topology $\tau$,
 we usually consider the strong topology of $E$.
Also, a family of mappings $\{ T(p) : p \in \RealNumber_+^n \}$ is called
 an $n$-parameter $\tau$-continuous semigroup
 of nonexpansive mappings on $C$
 (in short, an $n$-parameter nonexpansive semigroup)
 if (sg~1), (sg~2) and the following (sg~3) are satisfied:
\begin{enumerate}
\setcounter{enumi}{2}
 \renewcommand{\labelenumi}{(sg \arabic{enumi})}
\item
 for each $p \in \RealNumber_+^n$,
 $T(p)$ is a nonexpansive mapping on $C$.
\end{enumerate}
We know that
 an $n$-parameter nonexpansive semigroup $\{ T(p) : p \in \RealNumber_+^n \}$
 has a common fixed point
 in the case that $E$ is uniformly convex and
 $C$ is bounded, closed and convex;
 see Browder \cite{REF:Browder1965_ProcNAS_3}.
Moreover, in 1974, Bruck \cite{REF:Bruck1974_Pacific} proved that
 an $n$-parameter nonexpansive semigroup $\{ T(p) : p \in \RealNumber_+^n \}$
 has a common fixed point
 in the case that $C$ is weakly compact, convex, and
 has the fixed point property for nonexpansive mappings.

Very recently, the author proved the following
 in \cite{REF:TSP_cfps1_1_04}.

\begin{thm}[\cite{REF:TSP_cfps1_1_04}]
\label{THM:1-para}
Let $E$ be a Banach space and let $\tau$ be a Hausdorff topology on $E$.
Let $\{ T(t) : t \geq 0 \}$ be
 a $1$-parameter $\tau$-continuous semigroup of mappings
 on a subset $C$ of $E$.
Let $\alpha$ and $\beta$ be positive real numbers
 satisfying $\alpha/\beta \notin \RationalNumber$.
Then
 $$ \bigcap_{t \geq 0} F \big( T(t) \big)
 = F \big( T(\alpha) \big) \cap F \big( T(\beta) \big) $$
 holds.
\end{thm}

\noindent
Using this theorem,
 for an $n$-parameter $\tau$-continuous semigroup
 $\{ T(p) : p \in \RealNumber_+^n\}$
 of mappings,
 we obtain
 $$ \bigcap_{p \in \RealNumber_+^n} F \big( T(p) \big)
 = \bigcap_{k=1}^n \Big( F \big( T(e_k) \big) \cap
  F \big( T(\sqrt{2} e_k) \big) \Big) . $$
That is,
 the set of common fixed points of $\{ T(p) : p \in \RealNumber_+^n \}$
 is the set of common fixed points of $2 n$ mappings.

In this paper, motivated by the above thing,
 we prove the direct generalization of Theorem \ref{THM:1-para}
 which says that
 the set of common fixed points of $\{ T(p) : p \in \RealNumber_+^n \}$
 is the set of common fixed points of some $n+1$ mappings.
To prove it, we use Kronecker's theorem (Theorem \ref{THM:Kronecker}).
We also discuss convergence theorems
 to a common fixed point of
 $n$-parameter nonexpansive semigroups $\{ T(p) : p \in \RealNumber_+^n \}$.

\section{Preliminaries}
\label{SC:preliminaries}

In this section, we give some preliminaries.
For a real number $t$,
 we denote by $[t]$
 the maximum integer not exceeding $t$.
It is obvious that $0 \leq t - [t] < 1$ for all $t \in \RealNumber$.

We use two kinds of the notions of linearly independent in this paper.
We recall that
 vectors $\{ p_1, p_2, \cdots, p_n \}$
 is linearly independent in the usual sense
 if and only if
 there exist no $(\lambda_1, \lambda_2, \cdots, \lambda_n) \in \RealNumber^n$
 such that
 $$ (\lambda_1, \lambda_2, \cdots, \lambda_n) \neq 0
 \quad\text{and}\quad
 \lambda_1 p_1 + \lambda_2 p_2 + \cdots + \lambda_n p_n = 0 . $$
On the other hand,
 we call that
 real numbers $\{ \alpha_1, \alpha_2, \cdots, \alpha_n \}$ is
 linearly independent over $\RationalNumber$
 if and only if
 there exist no $(\nu_1, \nu_2, \cdots, \nu_n) \in \Integer^n$
 such that
 $$ (\nu_1, \nu_2, \cdots, \nu_n) \neq 0
 \quad\text{and}\quad
 \nu_1 \alpha_1 + \nu_2 \alpha_2 + \cdots + \nu_n \alpha_n = 0 . $$
For example,
 $$ \left\{ 1, \sqrt{2}, \sqrt{3}, \sqrt{5}, \sqrt{7}, \sqrt{11},
 \sqrt{13}, \sqrt{17}, \sqrt{19}, \sqrt{23} \right\} $$
 is linearly independent over $\RationalNumber$.
For each irrational number $\gamma$, $\{ 1, \gamma \}$
 is also linearly independent over $\RationalNumber$.
The following theorem is Kronecker's theorem;
 see \cite{REF:Hardy_Wright} and others.

\begin{thm}[Kronecker, 1884]
\label{THM:Kronecker}
Let $\alpha_1, \alpha_2, \cdots, \alpha_n \in \RealNumber$ such that
 $\{ 1, \alpha_1, \alpha_2, \cdots, \allowbreak \alpha_n \}$ is
 linearly independent over $\RationalNumber$.
Then the set of cluster points of the sequence
 $$ \Big\{ \Big( k \alpha_1 - [k \alpha_1], \;
 k \alpha_2 - [k \alpha_2], \; \cdots, \;
 k \alpha_n - [k \alpha_n] \Big)
 : k \in \NaturalNumber \Big\} $$
 is $[0,1]^n$.
\end{thm}

Let $E$ be a Banach space.
We recall that
 $E$ is called strictly convex if
 $\| x+y \| / 2 < 1$
 for all $ x, y \in E $ with $ \| x \| = \| y \| = 1 $
 and $ x \neq y $.
$E$ is called uniformly convex if
 for each $\varepsilon > 0$,
 there exists $\delta > 0$ such that
 $\| x+y \| / 2 < 1 - \delta$
 for all $ x, y \in E $ with $ \| x \| = \| y \| = 1 $
 and $\| x - y \| \geq \varepsilon$.
It is obvious that a uniformly convex Banach space is strictly convex.
The norm of $E$ is called Fr\'echet differentiable if
 for each $x \in E$ with $\| x \| = 1$,
 $\lim_{t \rightarrow 0}
 (\| x + t y \| - \| x \|)/t$ exists and is attained
 uniformly in $y \in E$ with $\| y \| = 1$.
The following lemma is the corollary of
 Bruck's result in \cite{REF:Bruck1973_TransAMS}.

\begin{lem}[Bruck \cite{REF:Bruck1973_TransAMS}]
\label{LEM:Bruck}
Let $C$ be a subset of a strictly convex Banach space $E$.
Let $\{ T_1, T_2, \cdots, T_\ell \}$ be
 a family of nonexpansive mappings from $C$ into $E$
 with a common fixed point.
Let $\lambda_1, \lambda_2, \cdots, \lambda_\ell \in (0,1]$
 such that $\sum_{j=1}^\ell \lambda_j = 1$.
Then
 a mapping $S$ from $C$ into $E$ defined by
 $$ S x = \lambda_1 T_1 x + \lambda_2 T_2 x + \cdots + \lambda_\ell T_\ell x $$
 for $x \in C$ is nonexpansive and
 $$ F(S) = F(T_1) \cap F(T_2) \cap \cdots \cap F(T_\ell) $$
 holds.
\end{lem}

\section{Main Results}
\label{SC:results}

In this section, we prove our main results.

\begin{thm}
\label{THM:main}
Let $E$ be a Banach space and let $\tau$ be a Hausdorff topology on $E$.
Let $\{ T(p) : p \in \RealNumber_+^n \}$ be
 an $n$-parameter $\tau$-continuous semigroup of mappings
 on a subset $C$ of $E$.
Let $p_1, p_2, \cdots, p_n \in \RealNumber_+^n$ such that
 $\{ p_1, p_2, \cdots, p_n \}$ is linearly independent in the usual sense.
Let $\alpha_1, \alpha_2, \cdots, \alpha_n \in \RealNumber$ such that
 $\{ 1, \alpha_1, \alpha_2, \cdots, \alpha_n \}$
 is linearly independent over $\RationalNumber$, and
 $$ p_0 = \alpha_1 p_1 + \alpha_2 p_2 + \cdots + \alpha_n p_n
 \in \RealNumber_+^n . $$
Then
 $$ \bigcap_{p \in \RealNumber_+^n} F \big( T(p) \big)
 = F \big( T(p_0) \big) \cap F \big( T(p_1) \big) \cap \cdots \cap
 F \big( T(p_n) \big) $$
 holds.
\end{thm}

To prove it, we need some lemmas.
In the following lemmas and the proof of Theorem \ref{THM:main},
 we let
 $$ z \in F \big( T(p_0) \big) \cap F \big( T(p_1) \big) \cap \cdots \cap
 F \big( T(p_n) \big) . $$
That is,
 $$ T(p_0)z = T(p_1)z = \cdots = T(p_n)z = z . $$
Also, we put
 $$ \ell = \max\Big\{ \big[ | \alpha_j | \big] + 1 : 1 \leq j \leq n \Big\}
 \in \NaturalNumber,
 \quad
 \beta_k = \alpha_k + \ell > 0 $$
 for $k \in \NaturalNumber$ with $1 \leq k \leq n$, and
 $$ p_0' = \beta_1 p_1 + \beta_2 p_2 + \cdots + \beta_n p_n
 \in \RealNumber_+^n . $$

\begin{lem}
\label{LEM:p}
$T(p_0')z=z$ holds.
\end{lem}

\begin{proof}
Since $z$ is a common fixed point of
 $\{ T(p_k) : k \in \NaturalNumber, 0 \leq k \leq n \}$,
 we have
 \begin{align*}
 T(p_0')z
 &= T(\beta_1 p_1 + \beta_2 p_2 + \cdots + \beta_n p_n)z \\*
 &= T(p_0 + \ell p_1 + \ell p_2 + \cdots + \ell p_n)z \\
 &= T(p_0) \circ T(p_1)^\ell \circ T(p_2)^\ell \circ \cdots \circ
  T(p_n)^\ell z \\*
 &= z .
 \end{align*}
This completes the proof.
\end{proof}

\begin{lem}
\label{LEM:01}
For every $(\lambda_1, \lambda_2, \cdots, \lambda_n) \in [0,1)^n$,
 $$ T(\lambda_1 p_1 + \lambda_2 p_2 + \cdots + \lambda_n p_n)z = z$$
 holds.
\end{lem}

\begin{proof}
We first show
 $\{ 1, \beta_1, \beta_2, \cdots, \beta_n \}$
 is linearly independent over $\RationalNumber$.
Assume that
 $$ \nu_0 + \nu_1 \beta_1 + \nu_2 \beta_2 + \cdots + \nu_n \beta_n = 0 $$
 for some $(\nu_0, \nu_1, \nu_2, \cdots, \nu_n) \in \Integer^n$.
Then by the definition of $\beta_j$,
 we obtain
 $$ (\nu_0 + \nu_1 \ell + \nu_2 \ell + \cdots + \nu_n \ell)
 + \nu_1 \alpha_1 + \nu_2 \alpha_2 + \cdots + \nu_n \alpha_n = 0 . $$
Since $\nu_0 + \nu_1 \ell + \nu_2 \ell + \cdots + \nu_n \ell \in \Integer$
 and $\{ 1, \alpha_1, \alpha_2, \cdots, \alpha_n \}$
 is linearly independent over $\RationalNumber$,
 we have
 $$ \nu_0 + \nu_1 \ell + \nu_2 \ell + \cdots + \nu_n \ell
 = \nu_1 = \nu_2 = \cdots = \nu_n = 0 . $$
{}From this, we also have $\nu_0 = 0$.
Therefore
 $\{ 1, \beta_1, \beta_2, \cdots, \beta_n \}$
 is linearly independent over $\RationalNumber$.
So, by Kronecker's theorem (Theorem \ref{THM:Kronecker}),
 there exists a sequence $\{ \ell_k \}$ in $\NaturalNumber$ such that
 $ \ell_k < \ell_{k+1} $ for $k \in \NaturalNumber$ and
 $$ \lim_{k \rightarrow \infty} \ell_k \beta_j - [\ell_k \beta_j]
 = \lambda_j $$
 for all $j \in \NaturalNumber$ with $1 \leq j \leq n$.
We next show
 $$ T \left( \sum_{j=1}^n \big( \ell_k \beta_j - [\ell_k \beta_j] \big)
  p_j \right) z = z $$
 for all $k \in \NaturalNumber$.
We define $T(p_j)^0$ is the identity mapping on $C$.
For each $k \in \NaturalNumber$,
 we have
 \begin{align*}
 & T \left( \sum_{j=1}^n \big( \ell_k \beta_j - [\ell_k \beta_j] \big)
  p_j \right) z \\*
 &= T \left( \sum_{j=1}^n \big( \ell_k \beta_j - [\ell_k \beta_j] \big)
  p_j \right) \circ T(p_1)^{[\ell_k \beta_1]}
  \circ T(p_2)^{[\ell_k \beta_2]} \circ \cdots
  \circ T(p_n)^{[\ell_k \beta_n]} z \\
 &= T \left( \sum_{j=1}^n \ell_k \beta_j p_j \right) z \\
 &= T \left( \sum_{j=1}^n \beta_j p_j \right)^{\ell_k} z \\*
 &= T \left( p_0' \right)^{\ell_k} z = z
 \end{align*}
 by Lemma \ref{LEM:p}.
Since
 $$ \lim_{k \rightarrow \infty}
  \sum_{j=1}^n \big( \ell_k \beta_j - [\ell_k \beta_j] \big) p_j
 = \sum_{j=1}^n \lambda_j p_j, $$
 we obtain the desired result.
\end{proof}

\begin{lem}
\label{LEM:infinity}
For every $(\lambda_1, \lambda_2, \cdots, \lambda_n) \in [0,\infty)^n$,
 $$ T(\lambda_1 p_1 + \lambda_2 p_2 + \cdots + \lambda_n p_n)z = z$$
 holds.
\end{lem}

\begin{proof}
By Lemma \ref{LEM:01}, we have
 \begin{align*}
 & T \left( \sum_{j=1}^n \lambda_j p_j \right) z \\*
 &= T \left( \sum_{j=1}^n \big( \lambda_j - [\lambda_j] \big) p_j \right)
  \circ T(p_1)^{[\lambda_1]} \circ T(p_2)^{[\lambda_2]}
  \circ \cdots \circ T(p_n)^{[\lambda_n]} z \\
 &= T \left( \sum_{j=1}^n \big( \lambda_j - [\lambda_j] \big) p_j \right) z \\*
 &= z .
 \end{align*}
This completes the proof.
\end{proof}

\begin{proof}[Proof of Theorem \ref{THM:main}]
We fix $p \in \RealNumber_+^n$.
Since $\{ p_1, p_2, \cdots, p_n \}$ is linearly independent in the usual sense,
 there exists $(\lambda_1, \lambda_2, \cdots, \lambda_n) \in \RealNumber^n$
 such that
 $$ p = \lambda_1 p_1 + \lambda_2 p_2 + \cdots + \lambda_n p_n . $$
Put
 $$ m = \max\Big\{ \big[ | \lambda_j | \big] + 1 : 1 \leq j \leq n \Big\}
 \in \NaturalNumber . $$
We note that $\lambda_j + m > 0$ for all $j$.
By Lemma \ref{LEM:infinity}, we obtain
 \begin{align*}
 T(p)z
 &= T \left( \sum_{j=1}^n \lambda_j p_j \right) z \\*
 &= T \left( \sum_{j=1}^n \lambda_j p_j \right)
  \circ T(p_1)^m \circ T(p_2)^m \circ \cdots \circ T(p_n)^m z \\
 &= T \left( \sum_{j=1}^n \big( \lambda_j + m \big) p_j \right) z \\*
 &= z .
 \end{align*}
Since $p \in \RealNumber_+^n$ is arbitrary,
 we obtain the desired result.
\end{proof}

Theorem \ref{THM:main} is the direct generalization
 of Theorem \ref{THM:1-para}.
We give the proof of Theorem \ref{THM:1-para}
 by using Theorem \ref{THM:main}.

\begin{proof}[Proof of Theorem \ref{THM:1-para}]
Put $p_1 = \beta$.
Since $p_1 \neq 0$, $\{ p_1 \}$ is linearly independent
 in the usual sense.
Put $\alpha_1 = \alpha / \beta \in \RealNumber \setminus \RationalNumber$.
We note that
 $\{ 1, \alpha_1 \}$ is linearly independent over $\RationalNumber$.
Put $ p_0 = \alpha_1 p_1 $.
Then by Theorem \ref{THM:main}, we obtain
 $$ \bigcap_{t \geq 0} F \big( T(t) \big)
 = F \big( T(p_0) \big) \cap F \big( T(p_1) \big)
 = F \big( T(\alpha) \big) \cap F \big( T(\beta) \big) . $$
This completes the proof.
\end{proof}

As another direct consequence of Theorem \ref{THM:main},
 we obtain the following.

\begin{cor}
\label{COR:prime}
Let $E$ be a Banach space and let $\tau$ be a Hausdorff topology on $E$.
Let $\{ T(p) : p \in \RealNumber_+^n \}$ be
 an $n$-parameter $\tau$-continuous semigroup of mappings
 on a subset $C$ of $E$.
Put $ \alpha_k $ the square root of the $k$-th prime number
 for $k \in \NaturalNumber$ with $1 \leq k \leq n$, and
 $$ p_0 = \alpha_1 e_1 + \alpha_2 e_2 + \cdots + \alpha_n e_n
 \in \RealNumber_+^n . $$
Then
 $$ \bigcap_{p \in \RealNumber_+^n} F \big( T(p) \big)
 = F \big( T(p_0) \big) \cap F \big( T(e_1) \big) \cap F \big( T(e_2) \big)
  \cap \cdots \cap F \big( T(e_n) \big) $$
 holds.
\end{cor}

Using Lemma \ref{LEM:Bruck}, we obtain the following.

\begin{cor}
\label{COR:sc}
Let $E$ be a strictly convex Banach space and
 let $\tau$ be a Hausdorff topology on $E$.
Let $\{ T(p) : p \in \RealNumber_+^n \}$ be
 an $n$-parameter $\tau$-continuous semigroup of nonexpansive mappings
 on a subset $C$ of $E$ with a common fixed point.
Let $\{ \alpha_1, \alpha_2, \cdots, \alpha_n \}$
 and $\{ p_0, p_1, p_2, \cdots, p_n \}$ as in Theorem \ref{THM:main}.
Define a nonexpansive mapping $S$ from $C$ into $E$ by
 $$ Sx = \lambda_0 T(p_0)x + \lambda_1 T(p_1)x + \cdots + \lambda_n T(p_n)x $$
 for $x \in C$,
 where $\lambda_0, \lambda_1, \lambda_2, \cdots, \lambda_n \in (0,1)$
 with $\sum_{j=0}^n \lambda_j = 1$.
Then
 $$ \bigcap_{p \in \RealNumber_+^n} F \big( T(p) \big) = F(S) $$
 holds.
\end{cor}

\begin{cor}
\label{COR:uc}
Let $E$ be a uniformly convex Banach space and
 let $\tau$ be a Hausdorff topology on $E$.
Let $\{ T(p) : p \in \RealNumber_+^n \}$ be
 an $n$-parameter $\tau$-continuous semigroup of nonexpansive mappings
 on a bounded closed convex subset $C$ of $E$.
Let $\{ \alpha_1, \alpha_2, \cdots, \alpha_n \}$
 and $\{ p_0, p_1, p_2, \cdots, p_n \}$ as in Theorem \ref{THM:main}.
Define a nonexpansive mapping $S$ on $C$
 as Corollary \ref{COR:sc}.
Then
 $$ \bigcap_{p \in \RealNumber_+^n} F \big( T(p) \big) = F(S) $$
 holds.
\end{cor}

\section{Convergence Theorems}
\label{SC:convergence}

Using Theorem \ref{THM:main},
 we can prove many convergence theorems
 to a common fixed point of
 an $n$-parameter $\tau$-continuous semigroup
 $\{ T(p) : p \in \RealNumber_+^n \}$ of nonexpansive mappings.
In this section, we state some of them.
In the following theorems,
 we always
 let $E$, $\tau$, $C$, $\{ T(p) \}$, $\{ p_j \}$, $\{ \alpha_j \}$ and
 $\{ \lambda_j \}$ as follows:

\begin{itemize}
\item
Let $E$ be a Banach space and
 let $\tau$ be a Hausdorff topology on $E$.
Let $\{ T(p) : p \in \RealNumber_+^n \}$ be
 an $n$-parameter $\tau$-continuous semigroup of nonexpansive mappings
 on a bounded closed convex subset $C$ of $E$.
Let $p_1, p_2, \cdots, p_n \in \RealNumber_+^n$ such that
 $\{ p_1, p_2, \cdots, p_n \}$ is linearly independent in the usual sense.
Let $\alpha_1, \alpha_2, \cdots, \alpha_n \in \RealNumber$ such that
 $\{ 1, \alpha_1, \alpha_2, \cdots, \alpha_n \}$
 is linearly independent over $\RationalNumber$, and
 $ p_0 = \alpha_1 p_1 + \alpha_2 p_2 + \cdots + \alpha_n p_n
 \in \RealNumber_+^n $.
Let $\lambda_0, \lambda_1, \lambda_2, \cdots, \lambda_n \in (0,1)$
 such that $\sum_{j=0}^n \lambda_j = 1$.
\end{itemize}

By the results of Bruck \cite{REF:Bruck1979_Israel}
 and Reich \cite{REF:Reich1979_JMAA},
 we obtain the following;
 see also Baillon \cite{REF:Baillon1975}.

\begin{thm}
\label{THM:Bruck-Reich}
Assume that $E$ is uniformly convex and
 the norm of $E$ is Fr\'echet differentiable.
Define a nonexpansive mapping $S$ on $C$ by
 $$ Sx = \sum_{j=0}^n \lambda_j T(p_j) x $$
 for all $x \in C$.
Define two sequences $\{ x_k \}$ and $\{ y_k \}$ in $C$ by
 $$ x \in C,
 \quad
 x_{k}
 = \frac{Sx + S^2 x + S^3 x + \cdots + S^k x}{k} $$
 for $k \in \NaturalNumber$, and
 $$ y_1 \in C,
 \quad
 y_{k+1} = \frac{1}{2} S y_k + \frac{1}{2} y_k $$
 for $k \in \NaturalNumber$.
Then
 $\{ x_k \}$ and $\{ y_k \}$ converge weakly to a common fixed point of
 $\{ T(p) : p \in \RealNumber_+^n \}$.
\end{thm}

By the results of Browder \cite{REF:Browder1967_ARMA}
 and Wittmann \cite{REF:Wittmann1992_ArchMath},
 we obtain the following;
 see also Halpern \cite{REF:Halpern1967_BullAMS}.

\begin{thm}
\label{THM:Browder-Wittmann}
Assume that $E$ is a Hilbert space.
Define a nonexpansive mapping $S$ on $C$ as Theorem \ref{THM:Bruck-Reich}.
Let $\{ s_k \}$ and $\{ t_k \}$ be sequences in $(0,1)$ satisfying
 $$ \lim_{k \rightarrow \infty} s_k
 = \lim_{k \rightarrow \infty} t_k = 0, \quad
 \sum_{k=1}^\infty t_k = \infty, \quad\text{and}\quad
 \sum_{k=1}^\infty | \; t_{k+1} - t_k \; | < \infty . $$
Define two sequences $\{ x_k \}$ and $\{ y_k \}$ in $C$ as
 $$ x_k = (1-s_k) S x_k + s_k u $$
 for $k \in \NaturalNumber$, and
 $$ y_1 \in C,
 \quad
 y_{n+1}
 = (1-t_k) \; S y_k + t_k \; u $$
 for $k \in \NaturalNumber$.
Then
 $\{ x_k \}$ and $\{ y_k \}$ converges strongly to a common fixed point of
 $\{ T(p) : p \in \RealNumber_+^n \}$.
\end{thm}

By the result of Rod\'{e} \cite{REF:Rode1982_JMAA},
 we obtain the following.

\begin{thm}
\label{THM:ave}
Assume that $E$ is a Hilbert space.
Define a sequence $\{ x_k \}$ in $C$ by
 $$ x \in C
 \quad\text{and}\quad
 x_{k}
 = \frac{\sum\left\{ T\left( \sum_{j=0}^n \nu_j p_j \right) x
  : \nu_j \in \{ 1, 2, \cdots, k \} \right\}}{k^{n+1}} $$
 for $k \in \NaturalNumber$.
Then
 $\{ x_k \}$ converges weakly to a common fixed point of
 $\{ T(p) : p \in \RealNumber_+^n \}$.
\end{thm}

By the result of Ishikawa \cite{REF:Ishikawa1979_Pacific},
 we obtain the following.

\begin{thm}
\label{THM:Ishikawa}
Assume that $C$ is compact.
Define mappings $S_j$ on $C$ by
 $$ S_j x = \frac{1}{2} T(p_j) x + \frac{1}{2} x $$
 for all $x \in C$ and $j = 0, 1, 2, \cdots, n$.
Let $x_1 \in C$ and define a sequence $\{ x_k \}$ in $C$ by
 $$ x_{k+1} = \left[ \prod_{k_{n}=1}^k \left[ S_n \prod_{k_{n-1}=1}^{k_{n}}
 \left[ S_{n-1} \cdots \left[ S_2 \prod_{k_1=1}^{k_2} \left[ S_1
 \prod_{k_0=1}^{k_1} S_0 \right] \right] \cdots \right] \right] \right] x_1 $$
 for $n \in \NaturalNumber$.
Then $\{ x_k \}$ converges strongly to a common fixed point of
 $\{ T(p) : p \in \RealNumber_+^n \}$.
\end{thm}

\section{Counterexample}
\label{SC:counterexample}

In Corollary \ref{COR:sc},
 we assume that $\{ T(p) : p \in \RealNumber_+^n \}$ has a common fixed point.
The following example says this assumption is needed.

\begin{exmp}
\label{EX:2-para}
Put $E=C=\RealNumber$ and
 let $\tau$ be the usual topology on $E$.
Define a $2$-parameter $\tau$-continuous semigroup
 $\{ T(p) : p \in \RealNumber_+^2 \}$ of nonexpansive mappings on $C$ by
 $$ T(\lambda_1 e_1 + \lambda_2 e_2) x
 = x + \lambda_1 - \lambda_2 $$
 for $\lambda_1, \lambda_2 \in [0,\infty)$ and $x \in E$.
Define a nonexpansive mapping $S$ on $C$ by
 $$ Sx
 = \frac{\sqrt{2}+\sqrt{3}+1}{6} T(\sqrt{2}e_1 + \sqrt{3}e_2)x
 + \frac{3-\sqrt{2}}{6} T(e_1)x
 + \frac{2-\sqrt{3}}{6} T(e_2)x $$
 for $x \in C$.
Then
 $$ \bigcap_{p \in \RealNumber_+^n} F \big( T(p) \big)
 = \varnothing
 \subsetneqq C = F(S) $$
 holds.
\end{exmp}

\begin{proof}
Since $F \big( T(e_1) \big) = \varnothing$ and $Sx = x$ for all $x \in C$,
 we obtain the desired result.
\end{proof}

\begin{ack}
The author wishes to express his sincere thanks
 to Professor Shige\-ki Aki\-yama in Niigata University
 for giving the valuable suggestions concerning Kronecker's theorem.
\end{ack}

\end{document}